\documentclass[11pt,reqno]{amsart}
\usepackage{a4wide}
\usepackage{hyperref}
\usepackage{color}
\usepackage{cite}
\usepackage{amsfonts,amssymb}
\usepackage{mathrsfs}
\usepackage{amstext}

\definecolor{green}{rgb}{0,0.6,0}
\definecolor{blue}{rgb}{0,0,1}

\hypersetup{colorlinks,
linkcolor=blue,
filecolor=green,
citecolor=green}

\theoremstyle{plain}

\newtheorem*{Cor*}{Corollary}

\newtheorem*{Thm*}{Theorem}

\newtheorem*{Prop*}{Proposition}
\theoremstyle{definition}

\newtheorem*{Rmk*}{Remark}

\newtheorem*{Ex*}{Example}

\newtheorem*{Qu*}{Question}

\theoremstyle{remark}

\theoremstyle{definition}

\newcommand{\x}{\times}

\newcommand{\p}{\partial}
\newcommand{\om}{\omega}

\newcommand{\pf}{\longrightarrow}

\newcommand{\N}{{\mathbb{N}}}
\newcommand{\Z}{{\mathbb{Z}}}
\newcommand{\R}{{\mathbb{R}}}
\newcommand{\C}{{\mathbb{C}}}

\newcommand{\LLL}{\mathscr{L}}

\renewcommand{\H}{\mathrm{H}}

\newcommand{\Id}{\mathrm{Id}}

\newcommand{\ind}{\mathrm{ind}}
\newcommand{\Ham}{\mathrm{Ham}}
\newcommand{\Supp}{\mathrm{Supp}}

\newcommand{\CF}{\mathrm{CF}}
\newcommand{\HF}{\mathrm{HF}}

\newcommand{\SH}{\mathrm{SH}}

\newcommand{\Spec}{{\rm Spec}}

\newcommand{\Crit}{{\rm Crit}}

\renewcommand{\AA}{\mathcal{A}}

\newcommand{\MM}{\mathcal{M}}

\newcommand{\beq}{\begin{equation}}
\newcommand{\beqn}{\begin{equation}\nonumber}
\newcommand{\eeq}{\end{equation}}

\newcommand{\bea}{\begin{equation}\begin{aligned}}
\newcommand{\bean}{\begin{equation}\begin{aligned}\nonumber}
\newcommand{\eea}{\end{aligned}\end{equation}}

\newcommand{\one}
{{{\mathchoice \mathrm{ 1\mskip-4mu l} \mathrm{ 1\mskip-4mu l}
\mathrm{ 1\mskip-4.5mu l} \mathrm{ 1\mskip-5mu l}}}}

\numberwithin{equation}{section}
\numberwithin{figure}{section}

\begin{document}
\title[Symplectic homology and Leafwise intersection points]{Symplectic homology of displaceable Liouville domains and Leafwise intersection points}
\author{Jungsoo Kang}
\address{Department of Mathematical Sciences\\
     Seoul National University, Seoul, Korea}
\address{Mathematisches Institut,
Westf\"alische Wilhelms-Universit\"at M\"unster, M\"unster, Germany}
\email{jungsoo.kang@me.com}
\begin{abstract}
 In this note we prove that the symplectic homology of a Liouville domain $W$ displaceable in the symplectic completion vanishes. Nevertheless if the Euler characteristic of $(W,\p W)$ is odd, the filtered symplectic homologies of $W$ do not vanish and give rise to leafwise intersection points on the symplectic completion of $W$ for a perturbation displacing $W$ from itself. In contrast to the existing results we can find a leafwise intersection point for a given period but its energy varies by period instead. 
\end{abstract}

\maketitle
\setcounter{tocdepth}{1}
\section{Main result}

Let $(W,d\lambda)$ be a Liouville domain, i.e. a compact exact symplectic manifold with contact type boundary $(V,\lambda|_V\!\!=\!\alpha)$. A neighborhood of $V$ in $W$ can be trivialized by the Liouville flow as $(V\x(1-\epsilon,1],d(r\alpha))$. Then the symplectic completion of $(W,d\lambda)$ is defined by
$$
\widehat W=W\cup_{\p W} V\x[1,\infty), \quad\widehat\om=\left\{\begin{array}{ll} d\lambda &\quad\textrm{on}\quad W,\\[1ex] d(r\alpha)&\quad\textrm{on}\quad V\x[1,\infty).\end{array}\right.
$$
We denote by $\widehat\lambda$ a primitive 1-form of $\widehat\om$ which is $\lambda$ on $W$ and $r\alpha$ on $V\x[1,\infty)$.
The contact manifold $V\x\{r\}$ for $r\in(1-\epsilon,\infty)$ is foliated by the leaves of the characteristic foliation spanned by the {\em Reeb vector filed} $R$ characterized by $\alpha(R)=1$ and $i_R d\alpha=0$. We recall that a closed Reeb orbit is {\em nondegenerate} if the linearized Poincar\'e return map associated to the orbit has no eigenvalue equal to 1. A contact manifold is called nondegenerate if all closed Reeb orbits are nondegenerate. We denote by $\varphi^t_R$ the flow of $R$. A submanifold $V$ in $(\widehat W,\widehat\om)$ is  said to be {\em displaceable} if there exists a perturbation $F\in C^\infty_c(S^1\x\widehat W)$ such that the associated Hamiltonian diffeomorphism $\phi_F$ defined below displaces $V$ from itself, i.e. $\phi_F(V)\cap V=\emptyset$. The displacement energy of $V$ in $\widehat W$ is defined by
$$
e(V):=\inf\{||F||\,|\,F\in C^\infty_c(S^1\x\widehat W),\,\phi_F(V)\cap V=\emptyset\}.
$$
We set $e(V)=\infty$ for the infimum of the empty set. Here by $||\cdot||$, we mean the Hofer norm which will be defined below.

We call $x\in V\x\{r\}$ a {\em leafwise intersection point} of $F$ if $\phi_F(\varphi_R^\eta(x))=x$ for some {\em period} $\eta\in\R$. For history of the leafwise intersection problem we refer to \cite{AF}. For a  given perturbation $F$, we abbreviate
$$
\frak{f}:=\max \pi(\Supp F),\quad
\pi:\left\{\begin{array}{l} W\setminus V\x(1-\epsilon,1]\pf\{0\},\\[1ex]
V\x(1-\epsilon,\infty)\pf(1-\epsilon,\infty).
\end{array}\right.
$$
Here $\pi$ is the projection along $V$ on $V\x(1-\epsilon,\infty)$. Our first result is:\\[-1.5ex]

\noindent\textbf{Theorem A.}
{\em Let $(V,\alpha)$ be a nondegenerate contact type boundary of a Liouville domain $(W,d\lambda)$ with $c_1(W)|_{\pi_2(W)}=0$. Suppose that $V$ is displaceable in $\widehat W$ and the Euler characteristic of $(W,\p W)$ is odd. Then for any $F\in C^\infty_c(S^1\x\widehat W)$ displacing $V$, there exists a leafwise intersection point of $F$ with $\tau$-period for every $\tau\in\R$ on $V\x\{r_\tau\}$ for some $r_\tau\in(1,\frak{f}]$.}\\[-1.5ex]

A generic starshaped hypersurface in $(\C^n,\om_{std})$ meets the requirements of the theorem.  To the best of  our knowledge, most results about the leafwise intersection problem concern the existence of or the number of leafwise intersection points with arbitrary periods on a fixed energy hypersurface. However in this note we find a fixed period leafwise intersection point but its energy could vary by period instead.

The theorem is an immediate consequence of the fact that the filtered symplectic homology $\SH^{(-\infty,\tau)}(W)$ does not vanish for every  $\tau\in\R$. Nevertheless we have:\\[-1.5ex]

\noindent\textbf{Corollary A1.}
{\em If $(V,\alpha)$ is nondegenerate and displaceable in $(\widehat W,\widehat\om)$, 
$$
i_*:\SH^{(-\infty,\tau)}(W)\to\SH^{(-\infty,\tau+e(V))}(W)
$$ 
induced by a canonical inclusion is a zero map for any $\tau\in\R$. In particular, the full symplectic homology $\SH(W)$ vanishes.}\\[-1.5ex]

After finishing the writing of the present paper, we became aware that Kai Cieliebak and Alexandru Oancea had obtained the corollary with a similar proof in their unfinished paper \cite{CO08}. There is an alternative proof of the vanishing of $\SH(W)$ by Ritter \cite{Rit} which makes use of a vanishing result of Rabinowitz Floer homology \cite{CF,AF}, the long exact sequence involving symplectic (co)homology and Rabinowitz Floer homology \cite{CFO}, and the unit of symplectic cohomology. However it seems to the author that the vanishing of the map $i_*:\SH^{(-\infty,\tau)}(W)\to\SH^{(-\infty,\tau+e(V))}(W)$ does not follow from his method. We will exploit the vanishing of $i_*$ in a forthcoming paper \cite{FK14}. Moreover an equivariant perturbation method in \cite{FS14} allows us to prove a  vanishing result of $S^1$-equivariant symplectic homology in the same way as the proof of Corollary A1. This result is also can be proved using big theorems, see \cite{BO12}.\\[-1.5ex]

\noindent\textbf{Corollary A2.}
{\em If $(V,\alpha)$ is nondegenerate and displaceable in $(\widehat W,\widehat\om)$, 
$$
i^{S^1}_*:\SH^{S^1,(-\infty,\tau)}(W)\to\SH^{S^1,(-\infty,\tau+e(V))}(W)
$$
induced by a canonical inclusion is a zero map for any $\tau\in\R$. In particular, the $S^1$-equivariant symplectic homology $\SH^{S^1}(W)$ vanishes.}

\section{Proof of the results}

\subsection{Convention and Notations.}
\begin{itemize}
\item Let a time-dependent almost complex structure $J_t$, $t\in S^1$ on $W$ be compatible with $\om$ and preserve the contact hyperplane field $\ker\alpha\subset TV$. We extend this on $\widehat W$ so that $J_t$ is invariant under the $\R_+$-action and $J_tr\p_r=R$ and $J_tR=-r\p_r$.  We call such almost complex structures {\em admissible}.
\item The {\em Hamiltonian vector field} $X_F$ associated to a Hamiltonian function $F\in C^\infty(S^1\x\widehat W)$ is defined by $i_{X_F}\widehat\om=dF$.
\item $\phi_F$ is the time one flow of $X_F$ and called a {\em Hamiltonian diffeomorphism}.
\item We denote by  $\Ham_c(\widehat W,\widehat\om)$ the group of Hamiltonian diffeomorphisms on $(\widehat W,\widehat\om)$ generated by compactly supported Hamiltonian functions.
\end{itemize}

\subsection{Symplectic homology and a perturbation}\quad\\[-1.5ex]

We first briefly recall symplectic homology and refer readers to \cite{BO,Vit} and references therein for further details. Since we have assumed that $(V,\alpha)$ is nondegenerate, $\Spec(V,\alpha)$ the set of all periods of closed Reeb orbits on $(V,\alpha)$ is a discrete subset in $\R_+=(0,\infty)$. We define an {\em admissible Hamiltonian} $H_\tau:\widehat W\to\R$ to have the following properties:
\begin{itemize}
\item[(i)] On $W\setminus V\x(1-\epsilon,1]$, $H_\tau$ takes values in $(-\epsilon,0)$ and is a $C^2$-small Morse function;
\item[(ii)] $H_\tau(r,x)=h_\tau(r)$ for some strictly increasing function $h_\tau(r)$ on $V\x(1-\epsilon,\infty)$ and  $h_\tau(r)=\tau r-\tau$ on $V\x(1,\infty)$;
\item[(iii)] $\tau\notin\Spec(V,\alpha)$ and $h_\tau''(r)>0$ on $(1-\epsilon,1)$.
\end{itemize}
We denote by $\LLL_{\widehat W}$ the component of contractible loops in $\widehat W$. We consider the action functional $\AA_{H_\tau}:\LLL_{\widehat W}\to\R$ defined by
$$
\AA_{H_\tau}(v)=-\int_{S^1}v^*\widehat\lambda-\int_{S^1}H_\tau(v)dt.
$$
There are two types of critical points of $\AA_{H_\tau}$: 
\begin{itemize}
\item[1)] critical points of the Morse function $H_\tau|_{W\setminus V\x(1-\epsilon,1]}$; 
\item[2)] solutions of 
\beq\label{eq:nontirivial critical point eq}
\p_t v=-h'_\tau(\pi\circ v)R(v).
\eeq
\end{itemize}
The second type solutions correspond to nondegenerate closed Reeb orbits with periods $h'_\tau(\pi\circ v(0))\in(0,\tau)$. Thus they  are in $V\x(1-\epsilon,1]$ and transversally nondegenerate due to (iii), i.e. 
$$
\ker [d\varphi_{R}^{-h'_\tau(\pi\circ v(0))}(v(0))-\one_{T_{v(0)}\widehat W}]=\langle\p_t v(0)\rangle,
$$ 
see \cite[Lemma 3.3]{BO}. Since we are considering parametrized periodic solutions, if $\gamma$ solves \eqref{eq:nontirivial critical point eq}, $S^1$-family of $\gamma(\cdot+t)$ do as well. We abbreviate $S_\gamma=\bigcup_{t\in S_1}\gamma(\cdot+t)$.  So $\AA_{H_\tau}$ is obviously not Morse but we still can define Floer homology of $\AA_{H_\tau}$ by Morse-Bott homology method, see \cite{Fra,BO}. We first choose a Morse function $f$ and a metric $g$ on the  critical manifold $\Crit\AA_{H_\tau}$ and associate the following index to critical points.
$$
\mu:\Crit f\pf \Z,\quad
\left\{\begin{array}{ll}\mu(\gamma_q)=\mu_{CZ}(\gamma)+\ind_f(\gamma_q)\qquad&\gamma_q\in\Crit f\cap S_\gamma\\[1ex]
\mu(p)=\ind_{-H_\tau}(p)-\frac{\dim W}{2}  &\,\;p\in\Crit H_\tau\subset\Crit f
\end{array}\right.
$$
where $\mu_{CZ}$ and $\ind_f$ stand for the Conley-Zehnder index and the Morse index respectively. Here we use $c_1(W)|_{\pi_2(W)}=0$ to obtain the $\Z$-valued index function $\mu$. We define the Floer chain group $\CF_n(H_\tau,f)$ by the $\Z/2$-vector space generated by critical points of $f$ and $H_\tau$ with index $n\in\Z$. The boundary operator is defined by counting gradient flow lines of $f$ together with cascades (that is, gradient flow lines of $\AA_{H_\tau})$. To be specific, for $\gamma_q^-,\gamma_q^+\in\Crit f$ and $m\in\N$, a {\em flow line from $\gamma_q^-$ to $\gamma_q^+$ with $m$ cascades} 
$$
(\textbf{w},\textbf{t})=\big((w_i)_{1\leq i\leq m},(t_i)_{1\leq i\leq m-1}\big)
$$
consists of $w_i\in C^\infty(\R\x S^1,\widehat W)$ solving
\beq\label{eq:gradient flow line}
\p_sw_i+J_t(w_i)(\p_tw_i-X_{H_\tau}(w_i))=0,
\eeq
the negative gradient flow equation of $\AA_{H_\tau}$ with respect to the metric $m$ on $T\LLL_{\widehat W}$ defined by
$$
m(v)[\xi,\zeta]=\int_0^1\widehat\om_v(\xi,J_t(v)\zeta)dt,\quad \xi,\,\zeta\in T_v\widehat W,
$$
and positive real numbers $t_i\in\R_+$ such that 
$$
\lim_{s\to\infty} (w_1(-s),w_m(s))\in W^u(\gamma_q^-;f)\x W^s(\gamma_q^+;f),\quad
\lim_{s\to-\infty} w_{i+1}(s)=\varphi_f^{t_i}(\lim_{s\to\infty}w_i(s))
$$
for $i=1,\dots,m-1$. Here $W^u(\gamma_q^-;f)$ resp. $W^s(\gamma_q^+;f)$ is the unstable manifold resp. the stable manifold and $\varphi_f^t$ is the flow of $-\nabla_g f$. We denote by $\widehat\MM_m(\gamma_q^-,\gamma_q^+)$ the space of flow lines with $m$ cascades from $\gamma_q^-$ to $\gamma_q^+$. We divide out the $\R^m$-action on this moduli space defined by shifting the cascades in the $s$-variable. We abbreviate 
$$
\MM_m(\gamma_q^-,\gamma_q^+)=\widehat\MM_m(\gamma_q^-,\gamma_q^+)/\R^m,\quad \MM(\gamma_q^-,\gamma_q^+)=\bigcup_{m\in\N\cup\{0\}}\MM_m(\gamma_q^-,\gamma_q^+).
$$
This moduli spaces is a smooth manifold of dimension $\dim\MM(\gamma_q^-,\gamma_q^+)=\mu(\gamma_q^-)-\mu(\gamma_q^+)-1$ for a generic admissible $J_t$, see \cite[Section 3]{BO}. 
Due to a maximum principle (\cite[Lemma 1.8]{Vit}) and a central theorem of Floer yield (\cite{Sal}) that  $\MM(\gamma_q^-,\gamma_q^+)$ is a finite set if $\mu(\gamma_q^-)-\mu(\gamma_q^+)=1$ (see also \cite{BO} for details in the Morse-Bott situation) and a boundary operator defined by
$$
\p:\CF_n^{(-\infty,b)}(H_\tau,f)\to\CF_{n-1}^{(-\infty,b)}(H_\tau,f),\quad \gamma_q^-\mapsto\sum_{\gamma_q^+\in\Crit f}\#_2\MM(\gamma_q^-,\gamma_q^+)\cdot \gamma_q^+
$$
indeed satisfies $\p\circ \p=0$. Here by $\#_2$ we mean the parity of the set. Now we define filtered symplectic homology by
$$
\SH_n^{(-\infty,b)}(W):=\lim_{\tau\to\infty}\HF_n^{(-\infty,b)}(H_\tau,f)
$$
where 
$$
\HF_n^{(-\infty,b)}(H_\tau,f):=\H_n(\CF_\bullet^{(-\infty,b)}(H_\tau,f),\p).
$$
Here the direct limit $\lim_{\tau\to\infty}$ is well-defined due to the maximal principle again. The full symplectic homology is defined as the direct limit of the above filtered ones, 
$$
\SH_n(W):=\lim_{b\to\infty}\SH_n^{(-\infty,b)}(W)=\lim_{\tau\to\infty}\lim_{b\to\infty}\HF_n^{(-\infty,b)}(H_\tau,f).
$$
	
 In order to study leafwise intersection points, we shall perturb the above action functional by $F\in C^\infty_c(S^1\x\widehat W)$. Before defining the perturbed action functional we make the time supports of $H_\tau$ and $F$ disjoint. To be precise, we pick a smooth function $\varrho\in C^\infty(S^1,[0,\infty))$ such that $\int_0^1\varrho dt=1$ and $\Supp\varrho\in(0,1/2)$, and we define $H^\varrho_\tau\in C^\infty(S^1\x\widehat W)$ by $H^\varrho_\tau(t,x)=\varrho(t)H_\tau(x)$. We can also modify $F$ so that $F(t,\cdot)=0$ for $t\in(0,1/2)$ and the Hamiltonian diffeomorphism $\phi_F$ remains unchanged.
$$
\AA_{H^\varrho_\tau+F}(v)=-\int_{S^1}v^*\widehat\lambda-\int_{S^1} H^\varrho_\tau(t,v)dt-\int_{S^1} F(t,v)dt.
$$
A critical point of $\AA_{H^\varrho_\tau+F}$ is a closed orbit $v\in\LLL_{\widehat W}$ solving
\beq\label{eq:critical point equation for the perturbed action functional}
\p_tv=X_{H^\varrho_\tau}(t,v)+X_{F}(t,v).
\eeq
For $v\in\Crit\AA_{H^\varrho_\tau+F}$, if $v(0)\in V\x(1-\epsilon,\infty)$, it satisfies $\phi_F(\varphi_R^{-h'_\tau(\pi\circ v(0))}(v(0)))=v(0)$ and thus $v(0)$ is a $-h'_\tau(\pi\circ v(0))$-period leafwise intersection point, see for instance \cite[Proposition 2.4]{AF}. On the other hand if $v(0)\in W$, $\phi_{H_\tau}(v(0))\in W$ as well and thus $\phi_F(W)\cap W\neq\emptyset$.
In a similar way to the unperturbed case, we can define $\HF_n(H^\varrho_\tau+ F)$ the Floer homology of $\AA_{H^\varrho_\tau+ F}$ for a generic perturbation $F$. We remark that since $\AA_{H^\varrho_\tau+F}$ is Morse for a generic $F$ in $C^\infty$-topology, see for instance \cite[Theorem 2.14]{AF} , we do not need an auxiliary Morse function. The following invariance property will play a crucial role.
\beq\label{eq:invariance property}
\SH_n^{(-\infty,\tau)}(W)=\HF_n^{(-\infty,\tau)}(H_\tau)\cong\HF_n(H^\varrho_\tau+F).
\eeq
The proof of this, which makes use of continuation homomorphisms, is fairly standard in Floer theory. We choose a cut-off function $\sigma:\R\to[0,1]$ so that $\sigma=0$ for $s\leq0$ and $\sigma=1$ for $s\geq1$. Then we set 
$$
F_s(t,x):=\sigma(s)F(t,x)\in C^\infty(\R\x S^1\x\widehat W).
$$
Then we consider solutions of 
\beq\label{eq:gradient flow equation for continuation homomorphism}
\left\{\begin{array}{l}
w:\R\x S^1\pf \widehat W\\[1ex]
w_-\in\Crit\AA_{H^\varrho_\tau},\quad w_+\in\Crit\AA_{H^\varrho_\tau+F},\quad \mu_{CZ}(w_-)=\mu_{CZ}(w_+)\\[1ex]
\p_sw+J(w)(\p_tw-X_{H^\varrho_\tau}(t,w)-X_{F_s}(t,w))=0\\
\end{array}\right.
\eeq
Here $w_\pm=\lim_{s\to\infty}w(\pm s)$. We note that nonconstant critical points of $\AA_{H_\tau^\varrho}$ still come in $S^1$-families and that critical points of $\AA_{H_\tau^\varrho+F}$ are still nondegenerate for a generic $F$. Let $f$ be a Morse function and $g$ be a metric on $\Crit\AA_{H_\tau^\varrho}$ as before. For $\gamma_q\in\Crit f$ and $v\in\Crit\AA_{H_\tau^\varrho+F}$, we consider a moduli space $\widehat\MM_m(\gamma_q,v)$ which is composed of 
$$
(\textbf{w},\textbf{t})=\big((w_i)_{1\leq i\leq m},(t_i)_{1\leq i\leq m-1}\big)
$$
such that $w_i\in C^\infty(\R\x S^1,\widehat W)$ is a solution of \eqref{eq:gradient flow line} (with respect to $H_{\tau}^\varrho$ instead of $H_{\tau}$) for $1\leq i\leq m-1$ and $w_m\in C^\infty(\R\x S^1,\widehat W)$ is a solution of \eqref{eq:gradient flow equation for continuation homomorphism} and that 
$$
\lim_{s\to-\infty} w_1(s)\in W^u(\gamma_q;f),\quad
\lim_{s\to-\infty} w_{i+1}(s)=\varphi_f^{t_i}(\lim_{s\to\infty}w_i(s)),\;\;t_i\in\R_+,\quad \lim_{s\to\infty}w_m(s)=v.
$$
As before we have the $\R^{m-1}$-action on $(w_i)_{1\leq i\leq m-1}$ and denote by 
$$
\MM_m(\gamma_q,v):=\widehat\MM_m(\gamma_q,v)/\R^{m-1},\quad \MM(\gamma_q,v):=\bigcup_{m\in\N}\widehat\MM_m(\gamma_q,v).
$$
Outside $W\cup(V\x[1,\frak{f}])$ where $X_F$ vanishes, $w$ must remain within $W\cup(V\x[1,\frak{f}])$ by the maximum principle again. Moreover the energy of $w$ is bounded in terms of the asymptotic data and $F$,
\bea\label{eq:energy estimate}
E(w):=\int_{-\infty}^\infty||\p_sw||_mds\leq\AA_{H^\varrho_\tau}(w_-)-\AA_{H^\varrho_\tau+F}(w_+)+||F||_-
\eea
where 
$$
 ||F||_-:=-\int_{S^1}\min_{x\in\widehat W}F(t,x)dt.
$$
Due to a uniform $C^0$-bound and a uniform energy bound on solutions of \eqref{eq:gradient flow equation for continuation homomorphism}, $\MM(\gamma_q,v)$ is a finite set if $\mu(\gamma_q)=\mu_{CZ}(v)$ and 
\bean
\Phi^{(-\infty,b)}:\CF_n^{(-\infty,b)}(\AA_{H^\varrho_\tau})&\pf\CF_n^{(-\infty,b+||F||_-)}(\AA_{H^\varrho_\tau+F}).\\
\gamma_q &\longmapsto\sum_{v\in\Crit\AA_{H_\tau^\varrho+F}}\#_2\MM(\gamma_q,v)\cdot v.
\eea
is a chain map. In a similar manner we have 
$$
\Psi^{(-\infty,b)}:\CF_n^{(-\infty,b)}(\AA_{H^\varrho_\tau+F})\pf\CF_n^{(-\infty,b+||-F||_-)}(\AA_{H^\varrho_\tau}).
$$
Due to a homotopy of homotopy argument  (see \cite{Sal} for details)
\beq\label{eq:homotopy of homotopy}
\Psi_*^{(-\infty,b+||F||_-)}\circ\Phi_*^{(-\infty,b)}=\iota^{b+||F||}_{b*}:\HF_n^{(-\infty,b)}(\AA_{H^\varrho_{\tau}})\to\HF_n^{(-\infty,b+||F||)}(\AA_{H^\varrho_{\tau}})
\eeq
where $||F||:=||-F||_-+||F||_-\geq0$ and $\iota^{b+||F||}_b$ is a canonical inclusion from $\CF_n^{(-\infty,b)}(\AA_{H^\varrho_\tau})$ into $\CF_n^{(-\infty,b+||F||)}(\AA_{H^\varrho_\tau})$. Taking the limit $b\to\infty$ we deduce that $\HF_n(H^\varrho_\tau,f)$ and $\HF_n(H^\varrho_\tau+F)$ are isomorphic. Moreover there is a canonical isomorphism between $\HF_n(H^\varrho_\tau,f)$ and $\HF_n(H_\tau,f)=\SH_n^{(-\infty,\tau)}(W)$ and thus \eqref{eq:invariance property} follows.
 It is worth pointing out that $\SH_n^{(-\infty,\tau)}(W)$ and $\HF_n(H_\eta+F)$ are not necessarily isomorphic if $\tau\neq \eta$.

One can play a similar game with $(-H_\tau,f)$ and a coboundary operator 
$$
\delta:\CF_{(-b,\infty)}^n(-H_\tau,f)\to\CF_{(-b,\infty)}^{n+1}(-H_\tau,f),\quad \gamma_q^+\mapsto\sum_{\gamma_q^-\in\Crit f}\#_2\MM(\gamma_q^-,\gamma_q^+)\cdot\gamma_q^-.
$$
We are able to define the filtered symplectic cohomology 
$$
\SH^n_{(-b,\infty)}(W):=\H^n(\CF^\bullet_{(-b,\infty)}(-H_\tau,f),\delta)
$$
and it is isomorphic to the perturbed filtered symplectic cohomology
$$
\SH^n_{(-\tau,\infty)}(W)\cong\HF^n(-H_\tau+F):=\H^n(\CF^\bullet(-H_\tau+F),\delta).
$$
for a generic $F\in C^\infty_c(S^1\x\widehat W)$.
\subsection{Proof of Theorem A}\quad\\[-1.5ex]

Since the Euler characteristic $\chi(W,\p W)$ is odd, $\HF(H_\tau,f)$ is nonzero for every admissible $H_\tau$'s. Indeed, if there is no closed Reeb orbit with period less than $\tau$ on $V$, $\HF(H_\tau,f)=\H(W,\p W)\neq0$. Even though there exist closed Reeb orbits with periods less than $\tau$, contributions of each closed Reeb orbit to the Floer chain groups of odd degree and to those of even degree are the same. That is, $\chi(\Crit f|_{S_\gamma})=0$ and thus we have
$$
\chi(\SH^{(-\infty,\tau)}(W))=\chi(\CF(H_\tau,f))=\chi(W,\p W)\neq0
$$ 
and hence $\SH^{(-\infty,\tau)}(W)$ is nonzero. 

Due to \eqref{eq:invariance property}, $\HF(H^\varrho_\tau+F)$ is nonzero for a generic $F\in C^\infty_c(S^1\x\widehat W)$ displacing $W$ as well. Thus we have at least one closed solution $v$ of \eqref{eq:critical point equation for the perturbed action functional}. If $v(0)\in \Sigma\x(\frak{f},\infty)$, the solution is nothing but a closed Reeb orbit whose period is $\tau\notin\Spec(V,\alpha)$ which never exist. On the other hand if $v(0)\in W$, the solution gives rise to a self-intersection point of $W$ by $F$ which contradicts to our displacing assumption on $V$. We remark that $V$ is displaceable in $\widehat W$ if and only if the whole filling $W$ is displaceable in $\widehat W$, see \cite[Lemma 3.4]{FSvK12}.

Therefore $v(0)$ has to lie on $\Sigma\x\{r\}$ for $r\in(1,\frak{f}]$ and this is a leafwise intersection point of $F$ with period $-\tau\notin -\Spec(V,\alpha)$. This argument goes through for the symplectic cohomology $\SH^n_{(-\tau,\infty)}(W)$ and thus we additionally obtain leafwise intersection points of $F$ with periods in $\R_+\setminus \Spec(V,\alpha)$.

Since $\Spec(V,\alpha)$ is discrete and $V\x(1,\frak{f}]$ is precompact, we can find leafwise intersection points in $\Sigma\x(1,\frak{f}]$ with periods in $\{0\}\cup\pm\Spec(V,\alpha)$ as well, see for instance \cite[p.130]{Kan}. Moreover since we have proved the theorem  for $F$ generic in $C^\infty$-topology, the theorem remains true for any perturbation in $C^\infty_c(S^1\x\widehat W)$ displacing $W$ away from itself. \hfill $\square$

\subsection{Proof of Corollary A1}\quad\\[-1.5ex]

Let $F\in C^\infty_c(S^1\x\widehat W)$ be a perturbation displacing $V$ from itself.
As we pointed out in the previous proof, if $v\in\Crit\AA_{H^\varrho_\tau+F}$, $\pi\circ v(0)\in (1,\frak{f}]$ and 
$$
\p_tv=-\tau \varrho(t)R(v)+X_F(t,v).
$$
Thus there exists a constant $C_F>0$ satisfying
\bean
\AA_{H^\varrho_\tau+F}(v)&=-\int_{S^1}\pi\circ v(0)\alpha(-\tau \varrho(t)R(v)+X_F(t,v))dt-\int_{S^1}(H_\tau^\varrho+F)(t,v)dt\\
&=\tau\pi\circ v(0)-(\tau\pi\circ v(0)-\tau)-\int_{S^1}\big(\pi\circ v(0)\alpha(X_F(t,v))+F(t,v)\big)dt\\
&\geq{\tau}-C_F.
\eea
Therefore for any $b>0$, there exists $\tau(b)>0$ such that $\CF^{(-\infty,b+||F||_-)}(\AA_{H_\tau(b)+F})$ and hence $\HF^{(-\infty,b+||F||_-)}(\AA_{H_\tau(b)+F})$ vanish. Thus the following map from \eqref{eq:homotopy of homotopy} is a zero map.
$$
i_{b*}^{b+||F||}=\Psi_*^{(-\infty,b+||F||_-)}\circ\Phi_*^{(-\infty,b)}:\HF^{(-\infty,b)}(\AA_{H_{\tau(b)}})\to\HF^{(-\infty,b+||F||)}(\AA_{H_{\tau(b)}}).
$$
Since this holds for any displacing Hamiltonian function, taking $\tau(b)\to\infty$, the map
$$
i_{b*}^{b+e(V)}:\SH^{(-\infty,b)}(W)\to\SH^{(-\infty,b+e(V))}(W),\quad b\in\R
$$
is zero and the first assertion of the corollary is proved. Taking the limit $b\to\infty$, we obtain,
$$
0=\Psi_*\circ\Phi_*=\Id_{\SH(W)}
$$
and this proves the vanishing of $\SH(W)$.\hfill$\square$

\subsubsection*{Acknowledgments} {\em I thank Urs Frauenfelder and Felix Schlenk for sharing with me their preprint \cite{FS14}. This paper is an outcome of insightful discussions with Urs Frauenfelder. I also thank Kai Cieliebak and Alex Oancea for sending their preprint \cite{CO08} where they independently proved Corollary A1. Many thanks to Peter Albers and Universit\"at M\"unster  for their warm hospitality. This work is supported by the National Research  Foundation of Korea (NRF) grant No. 2013R1A1A2004879 funded by Korea government (MEST) and by the SFB 878-Groups, Geometry, and Actions.}

\end{document}